\documentclass[11pt,leqno]{amsart}
\usepackage{amssymb, amsmath, amscd}
\usepackage{pstricks,pst-node}

\makeatletter
\renewcommand{\@seccntformat}[1]{\bf\@nameuse{the#1}.\quad}
\renewcommand\section{\@startsection{section}{1}%
            \z@{.7\linespacing\@plus\linespacing}{.5\linespacing}%
            {\normalfont\bfseries \boldmath}}
\renewcommand\subsection{\@startsection{subsection}{2}%
            \z@{.5\linespacing\@plus.7\linespacing}{-.5em}%
            {\normalfont\bfseries \boldmath}}
\renewcommand\subsubsection{\@startsection{subsubsection}{3}%
            \z@{.3\linespacing\@plus.5\linespacing}{-.5em}%
            {\normalfont\bfseries \boldmath}}
\makeatother

\theoremstyle{plain}
\newtheorem*{thm}{Theorem}
\newtheorem*{prop}{Proposition}
\newtheorem*{lemma}{Lemma}
\newtheorem*{cor}{Corollary}
\newtheorem*{guess}{Conjecture}
\newtheorem*{guessA}{Conjecture A}
\newtheorem*{guessB}{Conjecture B}

\theoremstyle{definition}

\newtheorem*{question}{Question}

\newtheorem*{claim}{Claim}

\theoremstyle{remark}

\newtheorem*{rem}{Remark}

\numberwithin{equation}{subsection}

\newcounter{listequation}

\def\note#1{{\small\tt <<#1>>}}  
\def\note#1{}              

\newcommand{\NN}{{\Bbb N}}
\newcommand{\ZZ}{{\Bbb Z}}
\newcommand{\RR}{{\Bbb R}}
\newcommand{\CC}{{\Bbb C}}
\newcommand{\ad}{\operatorname{ad}}

\newcommand{\Lie}{\operatorname{Lie}}

\newcommand{\tr}{\operatorname{tr}}
\renewcommand{\det}{\operatorname{det}}
\newcommand{\g}{\frak g}
\renewcommand{\k}{\frak k}
\newcommand{\h}{\frak h}
\renewcommand{\a}{\frak a}

\newcommand{\m}{\frak m}
\newcommand{\n}{\frak n}
\newcommand{\p}{\frak p}
\newcommand{\q}{\frak q}
\newcommand{\gl}{\frak g\frak l}
\renewcommand{\sl}{\frak s\frak l}
\newcommand{\so}{\frak s\frak o}
\renewcommand{\sp}{\frak s\frak p}
\newcommand{\su}{\frak s\frak u}
\renewcommand{\u}{\frak u}
\newcommand{\s}{\frak s}
\renewcommand{\t}{\frak t}
\newcommand{\z}{\frak z}
\renewcommand{\phi}{\varphi}
\renewcommand{\epsilon}{\varepsilon}

\begin{document}


\title[Multiplication of polynomials on Hermitian symmetric spaces]{Multiplication of polynomials on \\ Hermitian symmetric spaces and \\
Littlewood-Richardson coefficients}


\author{William Graham}
\address{Department of Mathematics, The University of Georgia,
Athens, GA 30602-7403, USA}
\email{wag@math.uga.edu}
\author{Markus Hunziker}
\address{Department of Mathematics, Baylor University, Waco, TX 76798-7328, USA}
\email{Markus\underline{\ }Hunziker@baylor.edu}

\begin{abstract}
Let $K$ be a complex reductive algebraic group and $V$ a representation of $K$. Let $S$ denote the ring of polynomials on $V$. Assume that the action of $K$ on $S$ is multiplicity free. If $\lambda$ denotes the isomorphism class of an irreducible representation of $K$, let $\rho_\lambda: K \rightarrow GL(V_{\lambda})$ denote the corresponding irreducible representation and $S_\lambda$ the $\lambda$-isotypic component of $S$. Write $S_\lambda \cdot S_\mu$ for the subspace of $S$ spanned by products of $S_\lambda$ and $S_\mu$. If $V_\nu$ occurs as an irreducible constituent of $V_\lambda\otimes V_\mu$, is it true that $S_\nu\subseteq S_\lambda\cdot S_\mu$? In this paper, the authors investigate
this question for representations arising in the context of Hermitian symmetric pairs. It is shown that the answer is yes in some cases and, using an earlier
result of Ruitenburg, that in the remaining classical cases, the answer is yes provided that a conjecture of Stanley on the multiplication of Jack polynomials is true. It is also shown how the conjecture connects multiplication in the ring $S$ to the usual Littlewood-Richardson rule.
\end{abstract}

\maketitle

\section{Introduction}

\subsection{}\label{question}
Let $K$ be a connected complex reductive group and let $X$ be a complex affine
algebraic variety with a $K$-action. We will assume that $X$ is a multiplicity free space, i.e., every finite dimensional irreducible representation of $K$ appears with multiplicity at most one in the algebra $S=\CC[X]$ of regular functions on $X$. Thus
$$
ÊÊÊ S \simeqÊ \bigoplus_{\lambda \in \Lambda} V_\lambda,
$$
where $\Lambda$ is a subset of $\widehat{K}$, the set of isomorphism classes of
finite dimensional irreducible $K$-representations, and $V_\lambda$ is a
representative of the class $\lambda\in \widehat{K}$. If $\lambda\in\Lambda$,
let $S_\lambda\subseteq S$ denote the $\lambda$-isotypic component of $S$, i.e.,
$S_\lambda$ is a $K$-stable subspace of $S$ such that $S_\lambda\simeq V_\lambda$. If $\lambda,\mu \in \Lambda$, let
$S_\lambda \cdot S_\mu$ be the subspace of $S$ spanned by the products of
elements in $S_\lambda$ and $S_\mu$. The subspace $S_\lambda\cdot S_\mu$
is $K$-stable and hence is a sum of certain $S_\nu$'s. We are interested in
describing the set of all $\nu\in \Lambda$ such that
$S_\nu \subseteq S_\lambda \cdot S_\mu$.
Clearly, a necessary condition forÊ $S_\nu \subseteq S_\lambda \cdot S_\mu$ is
that $V_\nu \hookrightarrow V_\lambda\otimes V_\mu$.

\begin{question}
Is it true that for $\nu,\lambda,\mu \in \Lambda$,
$$
ÊÊ S_\nu \subseteq S_\lambda \cdot S_\mu \quad \mbox{if and onlyÊ if}\quad
ÊÊÊ V_\nu \hookrightarrow V_\lambda\otimes V_\mu?
$$
\end{question}

\noindent
It is easy to find examples for which the answer is negative.
Consider the natural action of $K = SL_2(\CC)$ on $X=\CC^2$,
which is a multiplicity free action. More precisely, if $d$ is a non-negative integer, then the space $S_{d}$ of homogenous polynomials on $X$ of degree $d$
is an irreducible representation of $K$. In particular, $\Lambda$ can be identified with the set of non-negative integers. Clearly, $S_1\cdot S_1=S_2$  
and $V_1\otimes V_1 \simeq V_0\oplus V_2$. This counterexample can be fixed by
replacing $K=SL_2(\CC)$ with $SL_2(\CC)\times \CC^\times$, where the multiplicative group $\CC^\times$ acts by
multiplication on $X=\CC^2$. Then again $S_1\cdot S_1= S_2$, but now $V_1\otimes V_1\simeq \tilde{V}_0\oplus V_2$ with
$\tilde{V}_0\not\simeq V_0$, since $\CC^\times$ acts on $\tilde{V}_0$ with
non-trivial weight. This example suggests that a natural class of 
multiplicity free actions for which to study the question above is given
by the class of irreducible linear multiplicity free actions with the property
that the image of $K$ in $GL(X)$ contains the center $\CC^{\times}$
of  $GL(X)$. (A complete list of all such multiplicity free actions can be found in \cite{Kac1980}.) 
In this article we consider the nice subclass of multiplicity free actions that arise in the context of Hermitian symmetric spaces.

\subsection{}
Let $G_{\RR}/K_{\RR}$ be an irreducible Hermitian symmetric space of non-compact type and let $\g=\k\oplus\p^+ \oplus\p^-$Ê be the usual decomposition
of $\g=\Lie(G_\RR)\otimes_{\RR} \CC$ (see Section~\ref{S:hermitian} for more details).
The complexification of the compact group $K_\RR$ is a complex reductive group $K$ with one-dimensional center $\CC^{\times}$, and the action
of $K$ on $X=\p^+$ is multiplicity free. Table~1 shows a complete list of
the multiplicity free actions that arise in this way.
The last column of the table shows the rank of the Hermitian symmetric
space $G_{\RR}/K_{\RR}$ which is (by definition) the real rank of $G_\RR$. 

\begin{table}[h] 
\begin{center}
\begin{tabular}{|l|l|l|l|}
\hline
$G_\RR$ & $K$ & $X=\p^+$ & $r$ \\[.5ex]
\hline
$SU(p,q)$ & $S(GL_p(\CC)\times GL_q(\CC))$ & $\CC^p\otimes(\CC^q)^*$ &
$\min\{p,q\}$ \\[.5ex]
$Sp(n,\RR)$ & $GL_n(\CC)$ & $S^2(\CC^n)$ & $n$ \\[.5ex]
$SO^*(2n)$ & $GL_n(\CC)$ & $\wedge^2(\CC^n)$ & $\lfloor n/2\rfloor$ \\[.5ex]
$SO_{0}(n,2)$ & $SO_n(\CC) \times \CC^\times$ & $ \CC^n$ & $2$ \\[.5ex]
$E\, {I\!I\!I}$ÊÊ&Ê$Spin_{10}(\CC) \cdot \CC^\times$  & $\CC^{16}$ (spin)
& $2$\\[.5ex]
$E\, V\!I\!I$& $E_{6}(\CC) \cdot \CC^\times$ & $\CC^{27}$ (min)& $3$\\
\hline
\end{tabular}
\end{center}
\vspace{.5pc}
\caption{Multiplicity free actions associated to Hermitian symmetric spaces}\label{T:Table1}
\end{table}

We will study two conjectures. 

\begin{guessA}
For the multiplicity free action associated to the Hermitian symmetric space
$G_{\RR}/K_{\RR}$, i.e., for the $K$-action on  $X=\p^{+}$,
$$
ÊÊÊ S_\nu \subseteq S_\lambda \cdot S_\mu \quad
ÊÊÊ \mbox{if and only if}\quad V_\nu\hookrightarrow  V_\lambda\otimes V_\mu,
$$
where  $\lambda,\mu,\nu\in \Lambda$.
\end{guessA}

The second conjecture is closely related to the first and connects the problem of describing the set 
$\{\nu\in \Lambda \mid S_{\nu}\subseteq S_{\lambda}\cdot S_{\mu}\}$
to the classical Littlewood-Richardson rule. To state the conjecture, we need some more notation. Schmid proved in \cite{Schmid1969} that if $G_\RR/K_\RR$ has rank $r$,Ê then the non-zero $K$-isotypic components of $S=\CC[\p^+]$ are naturally parametrized by
partitions of length at most $r$, {\it i.e.\/}, the set
$$
ÊÊÊ \Lambda = \{(\lambda_1,\ldots,\lambda_r)\in\ZZ^r \mid
ÊÊÊ \lambda_1\geq \lambda_2\geq \cdots \geq \lambda_r\geq 0\}.
$$
This set $\Lambda$ also parametrizes the irreducible (polynomial)
representations of $GL_r(\CC)$ in the usual way. For $\lambda\in \Lambda$,
let $F_\lambda$ be the irreducible $GL_r(\CC)$-representation
of lowest weight $-\lambda$. For $\lambda,\mu,\nu\in \Lambda$,
we denote by $c^\nu_{\lambda\mu}$ the multiplicity of $F_\nu$ in
$F_\lambda\otimes F_\mu$. The numbers $c^\nu_{\lambda\mu}$ are
known as the Littlewood-Richardson coefficients for $GL_r(\CC)$.

\begin{guessB}\label{LR} 
Let $r$ be the rank of the Hermitian symmetric space $G_\RR/K_\RR$ and identify $\Lambda$ with the set of partitions of length at most $r$. Then for the $K$-action on $X=\p^+$,
$$
ÊÊÊ S_\nu \subseteq S_\lambda \cdot S_\mu \quad
ÊÊÊ \mbox{if and only if}\quad
ÊÊÊ c^{\nu}_{\lambda\mu}\not =0,
$$
where $\lambda,\mu,\nu \in \Lambda$ and $c^\nu_{\lambda\mu}$ are Littlewood-Richardson coefficients for $GL_r(\CC)$.
\end{guessB}

For the multiplicity-free actions
which correspond to Hermitian symmetric spaces of tube type (see Section \ref{SS: tube type}),
this conjecture was made by Ruitenburg \cite{Ruitenburg1989}, who observed that it holds for
$X = \CC^p \otimes (\CC^p)^*$.

The validity of Conjecture~B is connected to multiplication of Jack polynomials $P_\lambda^{(\alpha)}$, a well-known family of symmetric functions. 
Here $\lambda$ is a partition  and $\alpha$ is a real parameter. 
If $\alpha=1$, then $P_\lambda^{(\alpha)}=s_\lambda$ is the Schur symmetric function. There is an expansion
$$
P_\lambda^{(\alpha)} P_{\mu}^{(\alpha)} =
\sum_\nu f_{\lambda \mu}^{\nu}(\alpha) P_\nu^{(\alpha)},
$$
where the $f_{\lambda \mu}^{\nu}(\alpha)$ are rational functions in
the parameter $\alpha$. Note that since $P_\lambda^{(1)}=s_\lambda$, the value
$f_{\lambda \mu}^{\nu}(1)$ is the Littlewood-Richardson coefficient $c^{\nu}_{\lambda\mu}$.
In \cite{Stanley1989}, Stanley conjectured that the $f_{\lambda \mu}^{\nu}(\alpha)$
have a certain positivity property, which we can formulate as follows.
Let $c_{\lambda}(\alpha)$ and $c'_{\lambda}(\alpha)$ 
be the non-zero polynomials in $\alpha$ with non-negative integer coefficients as
defined in \cite[Chapter V, \S10]{Macdonald}. Then Stanley's conjecture
is that $c_{\lambda}(\alpha) c_{\mu}(\alpha) c'_{\nu}(\alpha) f_{\lambda \mu}^{\nu}(\alpha)$ are polynomials in $\alpha$ with non-negative integer coefficients.  Except for the non-negativity, this conjecture is true, by work of Knop and Sahi \cite{KnopSahi1997} and Lapoint and Vinet
\cite{LapointeVinet1995}. The main result of this paper, which we will prove in Section~\ref{S:Jack}, can then be phrased as follows.

\begin{thm} 
Suppose that Stanley's conjecture is true.
Then Conjecture~B is true for all Hermitian symmetric spaces.
\end{thm}

We also prove (without assuming Stanley's conjecture) that Conjecture~B is true
if $G_{\RR} = SU(p,q)$ or if the rank of the Hermitian symmetric space is at most $2$
(see Sections 2 and 3).

What about the validity of Conjecture~A? In Section~\ref{S:tensor}, we show the equivalence of Conjecture~A and Conjecture~B
for all Hermitian symmetric spaces of classical type.
More precisely, for $\lambda,\mu,\nu\in  \Lambda$ we express the multiplicity of $V_{\nu}$ in $V_\lambda\otimes V_\mu$ in terms of Littlewood-Richardson coefficients 
(see Table~\ref{T:multiplicities})
and show that $V_\nu\hookrightarrow  V_\lambda\otimes V_\mu$ if and only if $c^{\nu}_{\lambda\mu}\not =0$. In the classical cases corresponding to $G_{\RR}=Sp(n,\RR)$ and $SO^{*}(2n)$, this last equivalence follows from Klyachko's saturation conjecture, which is now a theorem of Knutson and Tao~\cite{KnutsonTao1999} and also of Derksen and Weyman~\cite{DerksenWeyman1990}.

\subsection{Acknowledgements.}
There is some overlap between the topics of this paper (especially Section~\ref{S:Jack}) and earlier work of Ruitenburg~\cite{Ruitenburg1989}.
For the reader's convenience, we have tried to give a relatively self-contained account. 

The second named author would like to thank Nolan Wallach and Friedrich Knop for their helpful comments.

\section{A classical example}

\subsection{}\label{regular rep}
Let $G$ be a complex reductive algebraic group. By reductive
we mean linearly reductive, i.e., all representations are completely
reducible.Ê (In particular, we do not assume here that $G$ is connected.)
Let $R = \CC[G]$ denote the ring of regular functions on $G$.
Recall that $R$ carries a natural action of $G \times G$,
the regular representation $\rho_{reg}$, given by
$$
(\rho_{reg}(g_1, g_2)f) (x) = f(g_1^{-1} x g_2)
$$
for $f \in R$ and $g_1, g_2, x \in G$.Ê By the algebraic
Peter-Weyl theorem, as a $G \times G$-representation
$R$ decomposes as a direct sum 
$$R = \bigoplus_{\lambda \inÊ \widehat{G}} R_{\lambda},$$
where $R_{\lambda} \simeq V_{\lambda} \boxtimes V_{\lambda}^*$.

The following proposition is observed in \cite[\S 3.2]{AlexeevBrion2004}.
\begin{prop}
Let $G$ be a complex reductive algebraic group and let $R=\CC[G]$
be the ring of regular functions on $G$.
Then 
 $$
Ê R_{\nu}\subseteq R_\lambda\cdot R_\mu
  \quad \mbox{if and only if}\quad 
  V_\nu\hookrightarrow  V_\lambda\otimes V_\mu
$$
for all $\lambda,\mu,\nu \inÊ \widehat{G}$.
\end{prop}

\begin{proof}
If $R_{\nu}\subseteq R_{\lambda} \cdot R_{\nu}$, then the $G \times G$-representation $V_{\nu} \boxtimes V_{\nu}^*$ appears in 
$(V_{\lambda} \boxtimes V_{\lambda}^*) \otimes (V_{\mu} \boxtimes V_{\mu}^*)$,
which implies that the $G$-representation $V_{\nu}$ appears in  $V_{\lambda} \otimes V_{\mu}$.
We now prove the converse.
For $\lambda\inÊ \widehat{G}$ let $\chi_{\lambda} \in R$ denote the character
of $V_{\lambda}$; i.e., $\chi_{\lambda}(g) =
\tr \rho_{\lambda}(g)$.
Recall that since $R_{\lambda}$ is spanned by the matrix coefficients of
$\rho_{\lambda}$ (with respect to any basis of $V_{\lambda})$, we have that $\chi_{\lambda} \in R_{\lambda}$.
Moreover, in $R$, we have $\chi_{\lambda} \cdot \chi_{\mu}
= \sum c_{\lambda \mu}^{\nu} \chi_{\nu}$, where
$c_{\lambda \mu}^{\nu}$ is the multiplicity of $V_{\nu}$
in $V_{\lambda} \otimes V_{\mu}$.Ê Let $\pi_{\lambda}:
R \rightarrow R_{\lambda}$ denote the projection.Ê Since
$R_{\lambda} \cdot R_{\mu}$ is a $G \times G$-stable subspace
of $R$, and since $G \times G$ is reductive, we have
$R_{\nu} \subseteq R_{\lambda} \cdot R_{\mu}$ if and only
if $\pi_{\nu}(R_{\lambda} \cdot R_{\mu})$ is non-zero.Ê If
$V_{\nu}$ occurs in $V_{\lambda} \otimes V_{\mu}$ then
$c_{\lambda \mu}^{\nu} \neq 0$.Ê Hence
$\pi_{\nu}( \chi_{\lambda} \cdot \chi_{\mu}) \neq 0$, so
$\pi_{\nu}(R_{\lambda} \cdot R_{\mu}) \neq 0$, completing the proof.
\end{proof}

\noindent
{\it Remark on notation.} Note that the decomposition $R=\sum R_\lambda$
is a decomposition into isotypic components for the $G\times G$-action
as well as the $G$-action induced by left multiplication on $G$
(the latter action being not multiplicity free).
We recall that the irreducible representations of $G\times G$
are parametrized by pairs $(\lambda,\mu)\inÊ \widehat{G}\times
Ê\widehat{G}$.
With this notation we then have
$R_{(\lambda,\mu)}=0$ if $\mu\not=\lambda^*$ and
$R_{(\lambda,\lambda^*)}=R_\lambda$, where $\lambda^*$
denotes the isomorphism class of the $G$-representation dual to
$V_\lambda$.

\subsection{}\label{p=q}
We use the notation of \ref{regular rep} with $G=GL_p(\CC)$. (To simplify notation, in the following we will write $GL_{p}$ for $GL_{p}(\CC)$.)
The set $\widehat{GL_p}$ is parametrized by integer sequences
$\lambda=(\lambda_1,\ldots,\lambda_p)$ with
$\lambda_1\geq \ldots \geq \lambda_p$.
More precisely, if $\lambda=(\lambda_1,\ldots,\lambda_p)$
let $F_\lambda$ be the irreducible $GL_p$-representation with lowest weight
$-\lambda_1\epsilon_1-\lambda_2\epsilon_2-\cdots-\lambda_p\epsilon_p$.
By abuse of notation
we will identify elements of $\widehat{GL_p}$ with such sequences.
Let $M_p$ denote the space of complex $p\times p$-matrices
and let $S=\CC[M_p]$. The open embedding $GL_p\hookrightarrow M_p$
induces a $GL_p\times GL_p$-equivariant embedding of coordinate rings
$S\hookrightarrow R=S[\det^{-1}]$.
This embedding identifies $S$
as the subring
$$
ÊS=\CC[M_p] \simeq \sum_{\lambda\in \Lambda} R_\lambda \ \subseteq\ R=\CC[GL_p]\ ,
$$
where
$\Lambda = \{(\lambda_1,\ldots,\lambda_p)\in\ZZ^p \mid
ÊÊÊ \lambda_1\geq \lambda_2\geq \cdots \geq \lambda_p\geq 0\}$.
Using notation as in the remark at the end of \ref{regular rep},
the only non-vanishing isotypic components for the
$GL_p\times GL_p$-action on $S$ are of the form
$S_{(\lambda,\lambda^*)}=R_\lambda$ for $\lambda\in \widehat{GL_p}$
with $\lambda_p\geq 0$.
It follows that Conjecture~A and Conjecture~B for the $K=GL_p\times GL_p$-action
on $X=M_p$ are an immediate consequence of the results in \ref{regular rep}.

\subsection{}\label{p<q}
Next we consider the action of $K=GL_{p}\times GL_{q}$ on $X=M_{p,q}$, 
where $M_{p,q}$ denotes the space of $p\times q$-matrices.
Without loss of generality, we may assume $p\leq q$. 
Let $M_{p,q} \twoheadrightarrow M_{p}$ be the natural projection 
given by ``forgetting the last $(q-p)$-columns''. This projection
is $GL_p$-equivariant with respect to the actions
given by left multiplication and thus we obtain a $GL_p$-equivariant
embedding $\CC[M_p] \hookrightarrow \CC[M_{p,q}]$. If $x_{ij}$ 
are the canonical coordinate functions on $M_{p,q}$ then 
$\CC[M_{p,q}]=\CC[x_{ij}\mid 1\leq i\leq p, 1\leq j\leq q]$ 
and $\CC[M_{p}]$ is identified (via the embedding given above) with
the subring $\CC[x_{ij}\mid 1\leq i, j\leq p]$.  
Let $N^-_p\subseteq GL_p$ be the group of lower diagonal unipotent $p\times p$-matrices and $N_q \subseteq GL_q$ the group of upper diagonal unipotent
$q\times q$-matrices. Let $\CC[M_{p,q}]^{N^-_p\times N_q}$ be the ring of 
$N^-_p\times N_q$-invariants in $\CC[M_{p,q}]$.
It is a theorem of classical invariant theory due to Weyl (see, e.g., 
\cite{GoodmanWallach}) that 
$\CC[M_{p,q}]^{N^-_p\times N_q} = \CC[u_1,\ldots, u_p]$,
where $u_k(x)$ is the $k$-th principal minor of $x\in M_{p,q}$.
Note that $u_k$ is in fact in $\CC[M_{p}]$ and hence
$$
Ê \CC[M_p]^{N^-_p\times N_p} = \CC[M_{p,q}]^{N^-_p\times N_q}.
$$
Furthermore, the polynomials $u_1,\ldots,u_r$ are algebraically independent
and for each 
$\lambda \in \Lambda = \{(\lambda_1,\ldots,\lambda_p)\in\ZZ^p \mid
ÊÊÊ \lambda_1\geq \lambda_2\geq \cdots \geq \lambda_p\geq 0\}$
, the polynomial
$u_\lambda=u_1^{\lambda_1}u_2^{\lambda_2 - \lambda_3}\cdots u_r^{\lambda_r}$
is a weight vector of weight
$(-\lambda_1\epsilon_1-\lambda_2\epsilon_2-\cdots-\lambda_p\epsilon_p,
\lambda_1\epsilon_1+\lambda_2\epsilon_2+\cdots+\lambda_p\epsilon_p )$
for the group $GL_p\times GL_q$.
Thus as a $GL_p\times GL_q$-representation,
$$
Ê \CC[M_{p,q}] \simeq \bigoplus_{\lambda\in \Lambda}
Ê F_\lambda^{(p)}\boxtimes \left(F^{(q)}_\lambda\right)^*
$$
Here we use superscripts $p$ and $q$ to emphasize that $F^{(p)}_\lambda$
is a representation of $GL_p$ and $F^{(q)}_\lambda$ is representation of
$GL_q$.

Let $S=\CC[M_p]$ and $T=\CC[M_{p,q}]$.
If $\lambda\in \Lambda$, let $S_\lambda$ and $T_\lambda$ denote the 
$\lambda$-isotypic components of $S$ and $T$, respectively, i.e.,
$S_\lambda \simeq F_\lambda^{(p)}\boxtimes (F^{(p)}_\lambda)^*$
as a $GL_p\times GL_p$-representation and $T_\lambda \simeq F_\lambda^{(p)}\boxtimes (F^{(p)}_\lambda)^*$
as a $GL_p\times GL_q$-representation. We note that 
$S_\lambda \hookrightarrow T_\lambda$ via the 
embedding $S\hookrightarrow T$. This follows from the observation 
that $S_\lambda$ and $T_\lambda$ are also the $\lambda$-isotypic components
for the $GL_p$-action on $S$ and $T$, respectively.

Suppose now that the $GL_{p}\times GL_{q}$-representation 
$F_\nu^{(p)}\boxtimes (F^{(q)}_\nu)^*$ appears in
$(F_\lambda^{(p)}\boxtimes (F^{(q)}_\lambda)^*)\otimes
(F_\mu^{(p)}\boxtimes (F^{(q)}_\mu)^*)
$.
Then the $GL_{p}$-representation
$F_\nu^{(p)}$ appears in $F_\lambda^{(p)}\otimes F_\mu^{(p)}
$.
By the previous subsection, 
it follows that $S_\nu\subseteq S_\lambda\cdot S_\mu$.
In particular, $u_\nu \in S_\lambda\cdot S_\mu$. 
Since $u_\nu\in T_\nu$ (via the embedding $S\hookrightarrow T$) 
this implies that $T_\nu\subseteq T_\lambda\cdot T_\mu$. It follows that
Conjecture~A and Conjecture~B are also true for the $GL_p\times GL_q$-action
on $X=M_{p,q}$ for $q>p$.

\begin{rem}
The argument in the previous section is a special case of a ``reduction to 
tube type'' argument introduced by Wallach in \cite{Wallach1979}.
We will discuss this kind of argument for general Hermitian symmetric pairs
in more detail in Section~\ref{tube type}.
\end{rem}

\section{Polynomials on Hermitian symmetric spaces}\label{S:hermitian}

\subsection{}
We recall some well-known results and constructions related to
Hermitian symmetric spaces which we will use in the following sections.
Let $\Omega$ be an irreducible Hermitian symmetric space of non-compact type.
Distinguishing a point $o\in \Omega$, we have $\Omega \simeq G_\RR/K_\RR$, 
where $G_\RR$ is the connected group of biholomorphic transformations of $\Omega$
and $K_\RR$ is the stabilizer group of $o$. The group $G_\RR$ is a simple
non-compact Lie group and $K_\RR$ is a maximal compact subgroup of $G_\RR$.
(These groups coincide with the groups of Table \ref{T:Table1}, up to local isomorphism.)  
Let $\g$ and $\k$ denote the complexified Lie
algebras of $G_\RR$ and $K_\RR$, respectively, and
let $\g=\k\oplus \p$ be the Cartan decomposition.
There is an element $z$ in the center of $\k$ such that
$\k=\CC z \oplusÊ [\frak k,\frak k]$ with $\ad(z)$ having the eigenvalues
$0$ and $\pm 1$ on $\frak g$.
Define $\frak p_{\pm} = \{x\in \g \mid [z,x]=\pm x\}$.
Then $\p=\p^+\oplus\p^-$, $[\p^\pm,\p^\pm]=0$, $[\p^\pm,\p^\mp]=\k$
and $[\k,\p^\pm]=\p^\pm$. 
Harish-Chandra contructed a canonical $K_\RR$-equivariant open embedding
$\Omega \hookrightarrow \p^+$ that sends $o$ to the origin in $\p^+$; via this embedding will view $\Omega$ as an open subset of $\p^{+}$.

\subsection{}
Let $\h_{\RR}$ be a Cartan subalgebra of $\k_{\RR}$.  The complexification $\h$ of $\h_{\RR}$ is a Cartan subalgebra of $\k$ and of $\g$.
Let $\Delta \subseteq \h^*$ be the root system of
$(\g,\h)$. For $\alpha\in \Delta$, let $\g_\alpha$ denote the root subspace
of $\g$ corresponding to $\alpha$.
Define the set of compact roots as
$\Delta_c = \{ \alpha\in \Delta \mid \g_\alpha \subseteq \k\}$
and the set of non-compact roots as
$\Delta_n = \{ \alpha\in \Delta \mid \g_\alpha \subseteq \p\}$.
We may choose a system of positive roots $\Delta^+$ for $\Delta$
such that $\p^+$ is the sum of the root subspaces for the roots
in $\Delta_n^+=\Delta_n^{\ }\cap\Delta^+$. The set
$\Delta_c^+=\Delta_c^{\ } \cap \Delta^+$ is a system of
positive roots for $\Delta_c$. 
Following Harish-Chandra we define a maximal set $\{\gamma_1,\ldots,\gamma_r\}$
of strongly orthogonal roots in $\Delta_n^+$ inductively as follows.
Let $\gamma_1$ be the largest root in $\Delta_n^+$ (with respect to the usual
ordering on $\Delta$ induced by the choice of $\Delta^+$); for $1<i\leq r$,
let $\gamma_i$ be the largest root in $\Delta^+_n$ that is
orthogonal to $\gamma_1,\ldots,\gamma_{i-1}$.
Table~\ref{gammas} shows the $\gamma_i$'s explicitly for all cases.
\begin{table}[h]
\begin{center}
\begin{tabular}{| l | l| l | l|}
\hline $\g$ & $\k$ & $r$ & $\{\gamma_1,\ldots , \gamma_r\}$ \\[1ex]
\hline $\sl_{p+q}(\CC)$ & $\s(\gl_p(\CC)\oplus \gl_q(\CC))$ & $\min\{p,q\}$ & $\{e_1-e_{p+q},e_2-e_{p+q-1},\ldots \}$ \\[1ex]
$\sp_{2n}(\CC)$Ê & $\gl_n(\CC)$ & $n$ & $\{2e_1,\ldots, 2e_n\}$Ê \\[1ex]
$\so_{2n}(\CC)$Ê & $\gl_n(\CC)$ & $\lfloor n/2 \rfloor$ & $\{e_1+e_2, e_3+e_4,\ldots \}$ \\[1ex]
$\so_{n+2}(\CC)$ & $\so_n(\CC) \oplus \CC$ & $2$ & $\{e_1 + e_2, e_1-e_2\} $ \\[1ex]
$\mathfrak e_6(\CC)$ & $\so_{10}(\CC) \oplus \CC$Ê & $2$ &
$\{\begin{smallmatrix}
ÊÊÊÊÊÊÊ 1&2&3&2& 1\\
ÊÊÊÊÊÊÊ & &2&
\end{smallmatrix},
\begin{smallmatrix}
ÊÊÊÊÊÊÊ 1&1& 1& 1&1\\
ÊÊÊÊÊÊÊ & & 0&
\end{smallmatrix}\}$Ê \\[1ex]
$\mathfrak e_7(\CC)$ & $\mathfrak e_6(\CC) \oplus \CC$& $3$ &
$\{\begin{smallmatrix}
ÊÊÊÊÊÊÊ 2&3&4&3&2&1\\
ÊÊÊÊÊÊÊ & &2&
\end{smallmatrix},
\begin{smallmatrix}
ÊÊÊÊÊÊÊ 0&1& 2&2&2&1\\
ÊÊÊÊÊÊÊ & & 1&
\end{smallmatrix}, \alpha_7\} $ \\[1ex]
\hline
\end{tabular}
\end{center}
\vspace{.5pc}
\caption{Strongly orthogonal roots}\label{gammas}
\end{table}
For $1\leq i\leq r$, define  $\pi_{i}:=\sum_{j=1}^{i}\gamma_{j}$;
the $\pi_{i}$ are  $\Delta_{c}^{+}$-integral and dominant weights.
Let $\Lambda = \{(\lambda_1,\ldots,\lambda_r)\in\ZZ^r \mid
ÊÊÊ \lambda_1\geq \lambda_2\geq \cdots \geq \lambda_r\geq 0\}$ be the set of partitions of length at most $r$. If $\lambda=(\lambda_1,\ldots,\lambda_r)\in \Lambda$, then $\sum_{i=1}^{r}\lambda_{i}\gamma_{i}=\sum_{i=1}^{r}(\lambda_{i}-\lambda_{i+1})\pi_{i}$ is a  $\Delta_{c}^{+}$-integral and dominant weight.
(Here it is understood that $\lambda_{r+1}=0$.)
InÊ \cite{Schmid1969}, Schmid gave an explicit decomposition of $S=\CC[\p^+]$
as a $K$-representation as follows.
If $\lambda\in \Lambda$, let $V_\lambda$ be the irreducible $\k$-representation
with lowest weight $-\sum_{i=1}^{r}\lambda_{i}\gamma_{i}$. Then
$$
ÊÊÊ S \simeqÊ \bigoplus_{\lambda \in \Lambda} V_\lambda,
$$
Furthermore, $S_\lambda$ is contained in the space of
homogeneous polynomials of degree $|\lambda|$, where $|\lambda|=\sum_{i=1}^r\lambda_i$.

\subsection{}\label{SS: tube type}
For the convenience of the reader and to provide a context for the material in Sections~\ref{tube type} and \ref{S:Jack}, we include here
some more details from Schmid~\cite{Schmid1969} and Koranyi-Wolf~\cite{KoranyiWolf1965}.
For each $\alpha\in \Delta_{n}^{+}$Ê let $\{h_\alpha, e_\alpha, e_{-\alpha}\}$ be an $\sl_2$-triple
such that $e_\alpha\in \g_\alpha$, $e_{-\alpha}\in \g_{-\alpha}$,
$h_\alpha\in [\g_\alpha,\g_\alpha]$ and $e_{-\alpha}=\sigma({e}_\alpha)$,
where $\sigma:\g\rightarrow \g$ denotes
complex conjugation with respect to $\g_\RR=\operatorname{Lie(G_{\RR})}$.
Let $e_+=e_{\gamma_1} + \ldots + e_{\gamma_r}$. Then $e_+$ is on the Shilov boundary of $\Omega$ in $\p^+$. In fact, the Shilov boundary is the $K_\RR$-orbit through $e_+$, and hence isomorphic to $K_\RR/M_{\RR}$, where $M_\RR=\operatorname{Stab}_{K_{\RR}}(e_{+})=\{k\in K_{\RR}\mid k e_{+}=e_{+}\}$.
Let $M=\operatorname{Stab}_{K}(e_{+})=\{k\in K \mid k e_{+}=e_{+}\}$
and let $\m$ be the Lie algebra of $M$.
Define an automorphism $\tau:\g\rightarrow \g$ by
$$
Ê \tau = \operatorname{Ad} \exp(i\pi (e_+ + \overline{e}_+)/2).
$$
Then $\tau^4=1$, $\tau^2(\k)=\k$, $\m=\{x\in \k\mid \tau(x)=x\}$
and the following are equivalent (see \cite[Proposition~4.4 and Lemma~4.8]{KoranyiWolf1965}):
\begin{itemize}
\item[(i)] $\Omega$ is a tube domain;
\item[(ii)] $\tau^2=1$;
\item[(iii)] $\tau(\k)=\k$;
\item[(iv)] $\dim \k -\dim \m = \dim \p^+$.
\end{itemize}

\noindent
Suppose from now on that $\Omega$ is a tube domain.
Then by (iv), the $K$-orbit through $e_+$ is open and dense in $\p^+$
and the corresponding open embeddingÊ $K/M \hookrightarrow \p^+$
gives a $K$-equivariant inclusion of coordinate rings
$$
ÊÊÊ S=\CC[\p^+] \hookrightarrow \CC[K/M].
$$
By the algebraic Peter-Weyl theorem
$$
ÊÊÊ \CC[K/M] \simeq \bigoplus_{\lambda \in \widehat{K}} V_\lambda\otimes (V_{\lambda}^*)^M.
$$
We may identify $\widehat{K}$ with the set of $\Delta^+_c$-dominant integral weights in $\h^*$. If $\lambda \in \widehat{K}$,
let $V_\lambda$ be the irreducible $K$-representation with lowest
weight $-\lambda$.
Then $V_\lambda^M\not =0$ (or equivalently $(V_{\lambda}^*)^M\not=0$)
if and only if $\lambda$
is in the lattice generated by the $\gamma_i$, i.e., if and only if
$\lambda$ is of the form $\lambda = \sum \lambda_i\gamma_i$ with
$\lambda_i\in \ZZ$ and $\lambda_1\geq \lambda_2 \geq \cdots \geq \lambda_r$
(see \cite[Lemma~1]{Schmid1969} ).
Furthermore, if $V_\lambda^M\not =0$ then $\dim V_\lambda^M=1$.
Schmid then showed (see \cite[Behauptung~c]{Schmid1969} )
that via the embedding $S=\CC[\p^+] \hookrightarrow \CC[K/M]$,
$$
S=\CC[\p^+]\simeq \bigoplus_{\lambda \in \Lambda} V_\lambda\otimes (V_\lambda^*)^M \hookrightarrow \CC[K/M],
$$
where $\Lambda = \{ \lambda \in \widehat{K} \mid
ÊÊÊ \lambda = \sum \lambda_i\gamma_i,\ \lambda_i\in \ZZ,\
ÊÊÊ \lambda_1\geq \lambda_2 \geq \cdots \geq \lambda_r\geq0\}$.
The connection with root systems of type $A$ is as follows. Since we assume that
$\Omega$ is of tube type, the pair $(\k,\m)$ is a symmetric pair with 
involution $\tau$. Let $\k_{\RR}=\m_{\RR}\oplus\s_{\RR}$ be the corresponding decomposition, i.e.,
$\m=\{x\in \k \mid \tau(x)=x\}$ and $\s= \{x\in \k \mid \tau(x)=-x\}$.  This complexifies to
$\k = \m \oplus \s$.
Define $\h_{\RR}^+=\h_{\RR} \cap \m_{\RR}$ and $\h_{\RR}^-=\h_{\RR} \cap \s_{\RR}$. Then $\h=\h^+\oplus \h^-$
and $\h^-$ is a maximal abelian
subspace of $\s$. By a result of Moore~\cite{Moore1964} (see also the remarks in the next section), if $\alpha\in \Delta_c^+$ then
$\alpha|_{\h^-} =0$ or $\alpha|_{\h^-} =\frac{1}{2}(\gamma_i-\gamma_j)$
with $i < j$. Thus the restricted root system $\Sigma$ is of type $A_{r-1}$.
Schmid's result above (which is a version of the Cartan-Helgason theorem; see also Section~\ref{S:Jack}) implies
that the $K$-representations occurring in $\CC[K/M]$ are exactly the 
$V_{\lambda}$, where $\lambda = \sum_{i=1}^{r} \lambda_{i}\gamma_i$ with $\lambda_i\in \ZZ$ and $\lambda_1\geq \lambda_2 \geq \cdots \geq \lambda_r$.
Table~\ref{T:KM} shows the symmetric pairs $(\k,\m)$ associated to the Hermitian symmetric spaces $\Omega=G_{\RR}/K_{\RR}$ of tube type. The last column shows the root multiplicity $m$ of the restricted roots. This root multiplicity will play 
an important role later in Section~\ref{S:Jack}.

\begin{table}
\begin{center}
\begin{tabular}{|l|l|l|l|}
\hline
$G_\RR$ & $\k$ & $\m$ & $m$ \\[1ex]
\hline
$SU(p,p)$ & $\frak{s}(\gl_p(\CC)\oplus \gl_p(\CC))$ & $ \sl_p(\CC) $  & $2$ \\[1ex]
$Sp(n,\RR)$ & $\gl_n(\CC)$ & ${\so}_n(\CC)$ & $1$ \\[1ex]
$SO^*(4p)$ & $\gl_{2p}(\CC)$ & $\sp_{\,2p}(\CC) $ & $4$ \\[1ex]
$SO_{0}(n,2)$, $n\geq 3$ & $\so_n(\CC)\oplusÊ\CC$& $\so_{n-1}(\CC)$Ê & $n-2$ \\[1ex]
$E\, V\!I\!I$ÊÊ & ${\frak e}_{6}(\CC) \oplus \CC$ & ${\frak f}_{4}(\CC)$ & $8$ \\[1ex]
\hline
\end{tabular}
\end{center}
\vspace{.5pc}
\caption{The symmetric pairs $(\k,\m)$ associated to 
$\Omega=G_{\RR}/K_{\RR}$ of tube type}\label{T:KM}
\end{table}

\begin{rem}
Note that in our context here, $K$ is not simply connected and 
$M$ is in general not connected. Furthermore, $K$ is not semisimple: 
it has a one-dimensional center. In Section~\ref{S:Jack}, to avoid
technical difficulties, we will work with symmetric pairs $(K,M)$, 
where $K$ is a connected and 
simply connected semisimple complex algebraic group. In \ref{SS:Main Theorem}, we will then return to the context of this section.
\end{rem}

\section{Reduction to tube type}\label{tube type}

\subsection{}
We retain the notation from the previous section. In particular,
let $\{\gamma_1, \ldots,\gamma_r\}$ be Harish-Chandra's strongly orthogonal
roots with the convention that $\gamma_1$ is the {\it largest} root in $\Delta_n^+$
and $\gamma_1>\gamma_2>\ldots>\gamma_r$. We point out to the reader that
this is different from the convention in much of the literature
(e.g., \cite{Moore1964}, \cite{Schmid1969}, \cite{Wallach1979}),
where $\gamma_1$ is taken to be the smallest root in $\Delta_n^+$.
If $(\g,\k)$ is of tube type (and only then), our $\gamma_i$Ê corresponds
to the other's $\gamma_{r-i}$.
Let $\h^-$ be the subspace of $\h$Ê spanned by the
coroots of $\gamma_1,\ldots,\gamma_r$ and let
$\h^+=\{h\in \h \mid \gamma_i(h)=0,\ 1\leq i \leq r\}$.
If $(\g,\h)$ is of tube type the spaces $\h^-$ and $\h^+$ agree with
the ones 
 \ref{SS: tube type}.
The following results are due to Moore~\cite{Moore1964}.
If $\alpha\in \Delta_c^+$ then
${\alpha|}_{\h^-}=\frac{1}{2}(\gamma_i-\gamma_j)$ with $i<j$,
${\alpha|}_{\h^-} = \frac{1}{2} \gamma_i$, or ${\alpha|}_{\h^-} = 0$.
If $\alpha\in \Delta_n^+$ then
$ \alpha|_{\h^-}=\frac{1}{2}(\gamma_i+\gamma_j)$ with $i\leq j$
or $ \alpha|_{\h^-} =Ê \frac{1}{2} \gamma_i$.
The Hermitian symmetric pair $(\g,\k)$ is of tube type if and only if
for every $\alpha\in\Delta$,
${\alpha|}_{\h^-}=\pm\frac{1}{2}(\gamma_i\pm\gamma_j)$, where
$1\leq i\leq j \leq r$.
The following construction is due to Wallach~\cite{Wallach1979}.
Let $\Delta_0=\{ \alpha \in \Delta \mid \ {\alpha|}_{\h^-}= \pm\frac{1}{2}(\gamma_i\pm\gamma_j)\}$.
Then $\Delta_0$ is a root subsystem of $\Delta$. Define $\tilde{\g}_0= \h\oplus \bigoplus_{\alpha\in \Delta_0} \g_\alpha$,
$\tilde{\k}_0=\k\cap\tilde{\g}_0$, and $\p_0^\pm=\p^\pm\cap \tilde{\g}_0$.
Then $\tilde{\g}_0=\tilde{\k}_0\oplus\p^+_0\oplus \p^-_0$,
$[\tilde{\k}_0,\p_0^\pm]=\p_0^\pm$,
$[\p_0^{\pm},\p_0^\pm]=0$ and
$[\p_0^\pm,\p_0^\mp]\subseteq \tilde{\k}_0$.
Define $\k_0=[\p^+_0,\p^-_0]$ and $\g_0=\k_0\oplus\p^+_0\oplus \p^-_0$.
By Lemma~2.2 in \cite{Wallach1979}, $\g_0$ is a simple
Lie subalgebra of $\g$ and $(\g_0, \k_0)$ is an irreducible Hermitian symmetric pair of tube type of rank $r$. Furthermore,
$\h^-\subseteq \k_0$ and $\h_0=(\h^+\cap \k_0)\oplus \h^-$ is a Cartan subalgebra
of $\k_0$ (and of $\g_0$).
By slight abuse of notation, the set
$\{\gamma_1,\ldots,\gamma_r\}$, via restriction to $\h_0$,
is also a maximal set of strongly orthogonal roots for $\p^+_0$.

\subsection{}
Let $\n^-=\bigoplus_{\alpha\in \Delta_c^+} \g_{-\alpha}$ and let
$\CC[\p^+]^{\n^-}$, i.e., the space of lowest vectors of the
$\k$-module.
By Schmid's result, the lowest weights in $\CC[\p^+]$ are of the form
$-\sum_{i=1}^r \lambda_i\gamma_i$. For $1\leq k\leq r$, let $u_k$
be a lowest weight vector in $\CC[\p^+]$ of weight
$-\sum_{i=1}^{k}\gamma_{i}$.
Then $\CC[\p^+]^{\n^-} =\CC[u_1,\ldots, u_r]$ and the
functions $u_1,\ldots,u_r$ are algebraically independent.
Recall that via the Killing form $(\p^+)^*\simeq \p^-$
and hence we can identify $\CC[\p^+]$ with the symmetric algebra $S(\p^-)$.
Similarly, we can identify $\CC[\p^+_0]$ with $S(\p^-_0)$.
Thus $\CC[\p^+_0]$ may be viewed as a subring of $\CC[\p^+]$.
Wallach showed (\cite{Wallach1979}, Lemma~3.3) that $u_k\in \CC[\p^+_0]$
for $1\leq k\leq r$. Using Schmid's result again, this time for the pair
$(\g_0,\k_0)$, we have
$$
Ê \CC[\p^+_0]^{\k_0\cap \n^-} = \CC[\p^+_0]^{\tilde{\k}_0\cap \n^-} = \CC[\p^+]^{\n^-}.
$$
Let $S=\CC[\p^+_0]$ and $T=\CC[\p^+]$.
If $\lambda=\sum_{i=1}^r \lambda_i\gamma_i$ then $S_\lambda = S\cap T_\lambda$.
Furthermore, $S_\lambda$ is an isotypic component of $S$ as a
$\k_0$-module as well as a $\tilde{\k}_0$-module.

\begin{prop}\label{reduction}
Let $(\g,\k)$ be an irreducible Hermitian symmetric pair of rank $r$
and let $(\g_0,\k_0)$ be the associated Hermitian symmetric pair of tube type of the same rank.
ThenÊ Conjecture~B is true for
$(\g,\k)$ if and only if it is true for $(\g_0,\k_0)$.
\end{prop}

\begin{proof}
We will show that $S_\nu\subseteq S_\lambda\cdot S_\mu$ if and only if
$T_\nu\subseteq T_\lambda\cdot T_\mu$.
By the remarks above, $S_\nu\subseteq S_\lambda\cdot S_\mu$ trivially
implies that $T_\nu\subseteq T_\lambda\cdot T_\mu$. To prove the converse we
will use an induction argument sketched by Enright and Wallach
in \cite{EnrightWallach2004}.
Define $\q=\tilde{\k}_0\oplus \u^+$, where $\u^+$ is the sum of
all root spaces $\g_\alpha$ with $\alpha\in \Delta_c^+$ such that
${\alpha|}_{\h^-}=\frac{1}{2}\gamma_i$ for some $1\leq i\leq r$.
Then $\q$ is a parabolic subalgebra of $\k$ with Levi factor
$\tilde{\k}_0$ and abelian nilradical $\u^+$ (see proof of Lemma~1
in \cite{EnrightWallach2004}). Let $\q^-=\tilde{\k}_0\oplus \u^-$
be the opposite parabolic.
If $E$ is an irreducible finite dimensional
$\tilde{\k}_0$-module, let $N(E)$
denote the $\k$-module obtained by
inducing from $\q^-$, i.e., $N(E)=U(\k)\otimes_{U(\q^-)} E$.

\begin{claim}
If $\nu$ is of the form $\nu=\lambda+\mu +\sum_{i=1}^r a_i\gamma_i$
with $\sum_{i=1}^r a_i =0$, then the inclusion
$S_\lambda\otimes S_\mu \hookrightarrow N(S_\lambda)\otimes N(S_\lambda)$ induces an equality
of $-\nu$-weight spaces
$$
Ê [S_\lambda\otimes S_\mu]_{-\nu} = [N(S_\lambda)\otimes N(S_\mu)]_{-\nu}
$$
\end{claim}

\noindent
{\it Proof of claim.\/} As an $\h$-module, $N(E)\simeq U(\u^+)\otimes E$.
The weights in $U(\u^+)$ restricted to $\h^-$ are of the form
$\sum_{i=1}^r \frac{1}{2}n_i\gamma_i$, where the $n_i$ are non-negative
integers. The restriction of the weights in $S_\lambda\otimes S_\mu$
to $\h^-$ are all of the form
$-\lambda-\mu + \sum_{i\leq j} \frac{1}{2}m_{ij} (\gamma_i-\gamma_j)$,
whereÊ the $m_{ij}$ are non-negative integers.

\medskip
\noindent
To show that $T_\nu\subseteq T_\lambda\cdot T_\mu$ implies
$S_\nu\subseteq S_\lambda\cdot S_\mu$ we take
$E=S_\lambda$ and $F=S_\mu$. By the universal property of
generalized Verma modules, there are canonical quotient maps
such $N(S_\lambda)\rightarrow T_\lambda$ and
$N(S_\mu)\rightarrow T_\mu$
such that the following diagram commutes:
$$
\CD
S_\lambda \otimes S_\mu @>>> S_\lambda\cdot S_\mu\\
@VVV @VVV\\
N(S_\lambda)\otimes N(S_\mu) @>>> T_\lambda\cdot T_\mu
\endCD
$$
If $\nu\in \Lambda$ such that $T_\nu\subseteq T_\lambda\cdot T_\mu$
then $|\nu|=|\lambda|+|\mu|$ and hence $\nu$ satisfies the hypothesis
of the claim. If we restrict the commutative diagram above
toÊ $-\nu$-weight spaces, the vertical arrow on the left is an
isomorphism. It follows that the inclusion
$[S_\lambda \cdot S_\mu]_{-\nu} \hookrightarrow [T_\lambda\cdot T_\mu]_{-\nu}$
is surjective and hence an isomorphism.
Since $S^{\tilde{\k}_0\cap \n^-}=T^{\n^-}$ we conclude that
$S_\lambda \cdot S_\mu$ contains a lowest weight vector of weight $-\nu$
and hence $S_\nu\subseteq S_\lambda\cdot S_\mu$.
\end{proof}

\subsection{}
Enright and Wallach~\cite{EnrightWallach2004} proved the following Pieri rule for multiplication of functions in $S=\CC[\p^+]$
by an induction on the rank of the Hermitian symmetric pair.

\begin{thm}[Enright-Wallach \cite{EnrightWallach2004}]\label{Pieri}
With notation as above, for every $\lambda\in \Lambda$ and $k\in \Bbb N$,
$$
ÊÊÊ S_\lambda \cdot S_{k\gamma_1} =ÊÊ \sum_\nu S_\nu,
$$
where the sum is over all $\mu\in \Lambda$ with $|\nu|=|\lambda|+k$ and
$\nu_1 \geq \lambda_1 \geq \nu_2 \geq \cdots \geq \lambda_{r-1} \geq \nu_r \geq \lambda_r .$
\end{thm}

\begin{cor}
Let $(\g,\k)$ be an irreducible Hermitian symmetric pair of rank $\leq 2$.
Then Conjecture~BÊ is true for $(\g,\k)$.
\end{cor}

\begin{proof}
By Proposition~\ref{reduction} we may assume that $(\g,\k)$ is of tube type.
In this case (see Schlichtkrull~\cite{Schlichtkrull1984}), the $\k$-module with lowest weight
$-\gamma_1-\cdots -\gamma_r$ is one-dimensional, i.e.,
$\dim S_{\gamma_1+\cdots+\gamma_r} =1$.
Since $S$ is a domain this implies that for every $\lambda \in \Lambda$,
$$
ÊÊÊ S_\lambda\cdot S_{\gamma_1+\cdots+\gamma_r} = S_{\lambda+\gamma_1+\cdots+\gamma_r}.
$$
Suppose now that $r=2$. Then for $\mu\in \Lambda$,
$$
ÊÊÊ S_{(\mu_1,\mu_2)} =
ÊÊÊ S_{(\mu_1-\mu_2,0)} \cdot S_{(\mu_2,\mu_2)}.
$$
By using the Pieri rule of the theorem above, we can then compute the decomposition of $S_\lambda\cdot S_\mu$ for any $\lambda,\mu \in \Lambda$.
If $\lambda,\mu,\nu\in\Lambda$ are such that $|\lambda|+|\mu|=|\nu|$, then
one finds that $S_{\nu}\subseteq S_\lambda\cdot S_\mu$ if and only $c^{\nu}_{\lambda\mu}\not=0$.
Explicitly, if  $a=\lambda_1-\lambda_2$, $b=\mu_1-\mu_2$ and $c=\nu_1-\nu_2$,
then
$$
c^\nu_{\lambda\mu}=
\begin{cases}
1 & \mbox{if $c=a+b-2k$ for some non-negative integer $k$}\\
0 & \mbox{otherwise.}
\end{cases}
$$
\end{proof}

\begin{rem}
In \cite{EnrightHunzikerWallach2004}, the theorem above was proved 
for $k=1$. This special case is enough to determine all the $K$-invariant 
ideals of $S=\CC[\p^+]$. The authors of \cite{EnrightHunzikerWallach2004}
were not aware at the time of writing that Ruitenburg in 
\cite{Ruitenburg1989} proved a Pieri rule and determined all the $K$-invariant
ideals of $S=\CC[\p^+]$ in a uniform manner for all cases corresponding
to Hermitian symmetric spaces of tube type. Ruitenburg did not work in 
the context of Hermitian symmetric spaces, but instead used the structure
of Riemannian symmetric spaces and spherical functions. In the next section
we will use the same approach.
\end{rem}

\section{Spherical functions and Jack polynomials}\label{S:Jack}

\subsection{} 
We begin with some preliminaries about symmetric spaces, and introduce some notation
which we will use throughout the section.
Let $K$ be a  connected and simply connected complex semisimple algebraic group.  Let $\tau$ be an involution of $K$ and let $M = K^{\tau}$; the group $M$ is connected
(see \cite[p. 171]{Loos}).   There exists a compact real form $(K_{\RR}, M_{\RR})$ of the pair
$(K,M)$; here $K_{\RR}$ is a maximal compact
subgroup of $K$ which is preserved by $\tau$, and $M_{\RR} = (K_{\RR})^{\tau}$
is a maximal compact subgroup of $M$.  This can be seen as follows.   There exists a real form $\k_0$ of $\k$ such that the involution $d \tau$ of $\k$ is the complexification of a Cartan involution
of $\k_0$ (see \cite[Lemma III.4.1]{HelgasonGGA}).  Let $\k_0 = \m_{\RR} + \s_0$ denote the corresponding Cartan decomposition; the complexification $\k = \m + \s$ is the decomposition of $\k$ into $+1$ and $-1$ eigenspaces for $\tau$.  Write $\s_{\RR} = i s_0$.  We can define a Cartan involution $d \tau'$ of $\k$ by requiring that $d \tau'$ act as multiplication by $1$ on $ {\frak k}_{\RR} :=  \m_{\RR} + \s_{\RR}$, and as multiplication by $-1$ on $i \k_0 +  i\s_{\RR}$.  Since $K$ is simply connected, there is a corresponding involution $\tau'$ of $K$, and the fixed point set  
$K_{\RR}=K^{\tau'}$ is a maximal compact subgroup of $K$ with Lie algebra ${\frak k}_{\RR}$ (see \cite[p. 252]{HelgasonDG}).  In particular, $K_{\RR}$ is connected and simply connected.  Moreover, by consideration of the Lie algebra one can see that $K_{\RR}$ is preserved by $\tau$.  Since $K_{\RR}$ is simply connected, the fixed point set $M_{\RR} := (K_{\RR})^{\tau'}$ is connected (again by \cite{Loos}).  Finally, $d\tau'$ restricts to a Cartan involution of $\m$, and $M_{\RR} = M^{\tau'}$, so $M_{\RR}$ is a maximal compact subgroup of $M$.  The space $K/M$ is a symmetric space, and it is
the complexification (in the sense of differential geometry) of the Riemannian symmetric space $K_{\RR}/M_{\RR}$.

\subsection{}\label{SS: spherical}
Let $\a_{\RR}$ be a maximal abelian subspace of $\s_{\RR}$, let $\h_{\RR}$ be any maximal abelian subspace of $\k_{\RR}$ containing $\a_{\RR}$, and let $\h$ denote the complex span of $\t_{\RR}$ in $\k$.  Then $\h$ is a Cartan subalgebra of $\k$ (\cite[p.~259]{HelgasonDG}).  Let $H \supseteq A$ denote the algebraic tori in $K$ whose Lie algebras are $\h$ and $\a$, respectively, and $H_{\RR} \supseteq A_{\RR}$ the compact tori of $K_{\RR}$ whose
Lie algebras are $\h_{\RR}$ and $\a_{\RR}$.  The torus $A$ is called a maximal split (or anisotropic) torus.  Let $X^*(A)$ denote the group of characters of $A$, viewed as a subset of $\a^*$; similarly we have $X^*(H) \subseteq \h^*$.   The group algebra over $\CC$ of $X^*(A)$ can be identified with the coordinate ring $\CC[A]$; write $e^{\lambda} \in \CC[A]$ for the element of $\CC[A]$ corresponding to the character $\lambda \in X^*(A)$.  

Let $\Sigma \subseteq X^*(A)$ denote the set of restricted roots; that is, the elements of $\Sigma$ are the non-zero weights for the $A$-action on $\k$.  Then $\Sigma$ is a root system in the real subspace of $\a^*$ it spans (see \cite[\S 4]{Richardson1982}).  Choose a positive system of roots $\Phi^{+}$ for 
$(\k,\h)$ and let $\Sigma^+$ denote the corresponding set of positive restricted roots.  This choice of positive system induces an ordering on $\a^*$ as usual, by the rule $\mu \leq \lambda$ iff $\lambda - \mu$ is a non-negative linear combination of positive restricted roots.  Let $X^*(H)^+\subseteq \h^{*}$ and $X^*(A)^+\subseteq \a^{*}$ denote the sets of dominant weights corresponding to the choice of positive system.  We write
$V_{\lambda}$ for the irreducible representation of $K$ with lowest weight $-\lambda \in X^*(H)^+$.
Note that the decomposition $\h = (\h \cap \m)\oplus \a$ allows us to view 
$\a^*$ as a summand in $\h^*$.  If  $\lambda \in 2X^*(A)$, then  $\lambda \in X^*(H)$; this follows since $A \cap M$ consists of the elements $a \in A$ with $a^2 = 1$.

\subsection{} Let $\CC[K]$ denote the coordinate ring of $K$.  By the algebraic Peter-Weyl theorem, 
$$
\CC[K/M] = \bigoplus_{\lambda \in \widehat{K}} V_{\lambda} \otimes (V_{\lambda}^*)^{M}.
$$
The group $M$ has the property that for any $\lambda \in \widehat{K}$, the dimension of the space $V^M_{\lambda}$ is either $1$ or $0$ (see 
\cite[Ch. IV \S 3]{HelgasonGGA}). 
The Cartan-Helgason theorem (see \cite[Ch.~V \S 2]{HelgasonGGA}, or \cite{Vust1974}) states that the irreducible $K$-representations occurring in $\CC[K/M]$ (that is, the $K$-representations with a non-zero $M$-fixed vector) are exactly the $V_{\lambda}$ where $\lambda \in 2 X^*(A)^+$.   
For each such weight $\lambda$ choose a left $M$-invariant function $\phi_{{\lambda}}$ in the $V_{{\lambda}}$-isotypic component of $\CC[K/M]$; this choice is unique up to scaling (we will choose a scaling in the discussion after Proposition \ref{p.jack}).  We may view $\phi_{{\lambda}}$ as a function on $K$, bi-invariant under $M$; these are called spherical functions.  The
spherical functions $\phi_{\lambda}$, as $\lambda$ runs through $2 X^*(A)^+$,
form a basis of the ring $\CC[K/M]^M$ as a vector space. Thus, we can expand the product of two spherical functions as follows:
$$
\phi_{{\lambda}} \phi_{{\mu}} =  \sum a_{\lambda \mu}^{\nu} \phi_{{\nu}},
$$
where the $a_{\lambda \mu}^{\nu}$ are constants. Let $R=\CC[K/M]$ and let $R_{\lambda}$ denote the $V_{\lambda}$-isotypic component of $R$. 
In \cite{Ruitenburg1989}, Ruitenburg proved the following result.
\begin{thm}[{Ruitenburg~\cite[Theorem~3.1]{Ruitenburg1989}}]
For fixed $\lambda,\mu,\nu \in 2 X^*(A)^+$,
$$
R_{\nu} \subseteq R_{\lambda} \cdot R_{\mu} \quad
\mbox{if and only if}\quad a_{\lambda \mu}^{\nu}\not=0.
$$ 
\end{thm}
To prove this result, Ruitenburg was using the compact real form described above. We will later use Ruitenburg' result to prove our main theorem from the introduction. 

\begin{rem}
If $\chi_{\lambda} \in \CC[K]$ denotes the character of the representation $\rho_{\lambda}: K \rightarrow GL(V_{\lambda})$, then up to scaling, $\phi_{\lambda}$ is the projection of $\chi_{\lambda}$ onto the space of $M$-invariants (with respect to the left action of $M$ on $\CC[K]$).  This follows from \cite[Theorem~4.2, Ch.~IV]{HelgasonGGA}. 
\end{rem}

\subsection{}
The negative of the Killing form induces a positive definite inner product on 
$\a_{\RR}$; let $s_{\alpha} \in GL(\a_{\RR})$ denote the reflection in the hyperplane $\alpha = 0$ in $\a_{\RR}$.
The Weyl group of the pair $(K,M)$, the ``little Weyl group'',  is $W = N_M(A)/Z_M(A) \simeq N_{M_{\RR}}(A_{\RR})/Z_{M_{\RR}}(A_{\RR})$. 
(Here, if $G \supseteq H$ are groups, 
$N_G(H)$ and $Z_G(H)$ denote the normalizer and centralizer of $H$ in $G$, respectively).  
The group $W$ acts on $\a_{\RR}$ and can be identified with its image in $GL(\a_{\RR})$, which is generated by the reflections $s_{\alpha}$.  (See \cite{Richardson1982} and \cite[Ch. VII]{HelgasonDG} for proofs of these facts.) 

The natural map $A/ (A \cap M) \rightarrow K/M$ induces a map of coordinate rings $\CC[K/M] \rightarrow \CC[A/(A \cap M)]$.  The induced map $\CC[K/M]^M \rightarrow \CC[A/(A \cap M)]^W$ is an isomorphism.  This is a reformulation of a theorem of Richardson {Richardson1982}.  Indeed, let $P$ be the subset of $K$ consisting of elements of the form $k \tau(k)^{-1}$.  There is a commutative diagram
$$
\begin{array}{ccc}
A / (A \cap M) & \to & K/M \\
\downarrow & & \downarrow \\
A & \to & P
\end{array}
$$
Here the horizontal maps are the inclusions.  The vertical maps are isomorphisms; the left vertical map takes $a (A \cap M)$ to $a^2$, and the right vertical map takes $kM$ to $k \tau(k)^{-1}$.  Note that $M$ acts by left translation on $K/M$ and by conjugation on $P$, and the right vertical map is $M$-equivariant.  Richardson's result states that the induced map $\CC[P]^M \to \CC[A]^W$ is an isomorphism, from which the version stated above follows.
Composing the isomorphisms $\CC[K/M]^M \rightarrow \CC[A/(A \cap M)]^W$ and
$\CC[A /(A \cap M)]^W \rightarrow \CC[A]^W$ yields an isomorphism $\CC[K/M]^M \to \CC[A]^W$.  
For $\lambda \in 2X^*(A)^+$, 
let $Q_{\lambda}$ denote the image of the spherical function 
$\phi_{{\lambda}}$. 
Finally, we define the function $P_{\lambda}$ in $\CC[A]^W$ by  $P_{\lambda}(a)=Q_{\lambda}(a^{-1})$ for $a\in A$. Thus, by definition,  we have
$$
  P_{\lambda}(a^2) =\phi_{\lambda}(a^{-1}).
$$
Note that in light of the remark at the end of the previous subsection, 
$P_{\lambda}=Q_{\lambda^{*}}$ (up to scaling), where $\lambda^{*}\in 2X^*(A)^+$
is such that $V_{\lambda^{*}}\simeq (V_{\lambda})^{*}$. 
Since the spherical functions $\phi_{\lambda}$ form a basis for $\CC[K/M]^M$,
the preceding discussion implies that the functions $P_{\lambda}$ 
form a basis for $\CC[A]^W$. Furthermore, the structure constants for the multiplication of the functions $P_{\lambda}$ and the spherical functions $\phi_{{\lambda}}$ are the same, i.e.,
$P_{{\lambda}} P_{{\mu}} =  \sum a_{\lambda \mu}^{\nu} P_{{\nu}}$, where
the $a_{\lambda \mu}^{\nu}$ are the same constants as at the end of the 
previous subsection.

\subsection{}\label{SS:p.jack}
In the cases most of interest in this paper, the functions $P_{\lambda}$ will turn out to be specializations of Jack polynomials.
Before we show this in the next subsection, we give an alternative characterization of the functions $P_{\lambda}$.
For $\lambda \in 2 X^*(A)^+$, define $m_{\lambda}\in \CC[A]^{W}$ by
$m_{\lambda}= (1/|W_{\lambda}|)\sum_{w \in W} e^{w \lambda / 2}$,
where $|W_{\lambda}|$ is the order of the stabilizer of $\lambda$ in $W$.
Define a function $\delta$ on the compact torus $A_{\RR}$ by $\delta(a) := \prod_{\alpha \in \Sigma^+} | 1- e^{\alpha}(a)|^{\operatorname{mult}(\alpha)}$, where $\operatorname{mult}(\alpha)$ is the multiplicity of the restricted root $\alpha$.
Finally, define an inner product on $\CC[A]$ by the rule
$$
\langle f, g \rangle_A = \int_{A_{\RR}}  f(a) \overline{g(a)} \delta(a)\, da.
$$
The following result is known, but because a complete proof seems hard to find in the literature, we provide the argument. As a side remark, it is perhaps interesting to note that most of the ideas  in the proof appear already in a paper of Harish-Chandra~\cite{Harish-Chandra} from 1958.

\begin{prop}\label{p.jack}
The functions $P_{\lambda}$ satisfy the following properties:
\begin{itemize}
\item[(a)] 
There exist constants $c_{\lambda\mu}$ such that 
$$ 
P_{\lambda} =  \sum_{\mu/2 \, \leq \, \lambda/2} c_{\lambda\mu} m_{\mu} 
\quad\text{with}\quad c_{\lambda\lambda}\not=0.
$$
\item[(b)] For all $\lambda\not=\mu$,
$$
\langle P_{\lambda},P_{\mu}\rangle_{A}=0.
$$  

\end{itemize}
\end{prop}

\begin{proof}
We prove the equivalent of the proposition for the functions $Q_{\lambda}$.
For $\lambda \in 2 X^*(A)^+$, define $n_{\lambda}\in \CC[A]^{W}$ by
$n_{\lambda}=  (1/|W_{\lambda}|) \sum_{w \in W} e^{-w \lambda / 2}$. 
The fact that $Q_{\lambda}=\sum c_{\lambda\mu} n_{\mu}$ with $c_{\lambda\lambda} 
\neq 0$ 
is proved in \cite[p.~274]{Harish-Chandra}.  Harish-Chandra \cite[p.~275]{Harish-Chandra} also gives a Freudenthal-type recursion formula for 
$\phi_{{\lambda}}$ which implies that
$Q_{\lambda} = \sum c_{\lambda\mu} n_{\mu}$, where the sum is over
$\mu\in 2X^{*}(A)^{+}$ such that $\mu/2 \in X^*(A)$ and $\mu/2 \leq \lambda/2$.  
(It is easy to see that the sum is over $\mu \leq  \lambda$, but the ordering is defined using integral linear combinations of positive roots, so the fact that we can divide this inequality by $2$ does not seem obvious without using Harish-Chandra's formula.  Arguments analogous to those in Humphreys~\cite[Lemma B, \S 13.3]{Humphreys}  show that the leading term of Harish-Chandra's formula is non-zero, so the formula determines $\phi_{{\lambda}}$ up to scaling.)  This implies (a).  As for (b), we have already observed that the $Q_{\lambda}$ form a basis of $\CC[A]^W$.  If $\lambda \neq \mu$, then
$$
\int_{K_{\RR}} \phi_{{\lambda}}(k) \overline{\phi_{{\mu}}(k)}\, dk = 0,
$$
by the usual orthogonality relation for matrix coefficients of different representations.
The function
$\phi_{{\lambda}} \overline{\phi_{{\mu}}}$ is bi-invariant under the group $M_{\RR}$,
so using the integration formula of Theorem 5.10 of \cite[Ch. I \S5]{HelgasonGGA}, we see that the above integral is
a non-zero constant times
$$
\int_{A_{\RR}/(A_{\RR} \cap M_{\RR})} \phi_{{\lambda}}(a) \overline{\phi_{{\mu}}(a)}  \delta(a^2) da.
$$
Via the isomorphism $A_{\RR}/(A_{\RR} \cap M_{\RR}) \to A_{\RR}$ which takes $a(A_{\RR} \cap M_{\RR})$ to $a^2$, this leads to the desired orthogonality.
\end{proof}

Henceforth, we will normalize the $\phi_{{\lambda}}$ so that the coefficient $c_{\lambda\lambda}$ equals $1$.

\subsection{} \label{SS:root}
Much of the following discussion generalizes to arbitrary root systems.  But for the application to our main results, we only need the case where the root system $\Sigma$ is of type $A_{r-1}$, so from now on, to simplify the exposition, we will assume this.  The torus $A$ is isomorphic to a maximal torus $T_{SL_{r-1}}$ via an isomorphism respecting the roots and the character lattices.  The reason is that our assumption that $K$ is simply connected implies that the integrality conditions characterizing the character lattice $X^*(A)$ are the same as the integrality conditions characterizing the character lattice $X^*(T_{SL_{r-1}})$ (see \cite[Ch.~VII \S 8]{HelgasonDG}).  Therefore, we can write the positive roots as
$\epsilon_{i}-\epsilon_{j}$, $1\leq i<j\leq r$, where $\epsilon_{i}$ are elements of $X^{*}(A)$ such that $\sum_{i=1}^{r} \epsilon_{i}=0$. 
Let $\gamma_{i}:=2\epsilon_{i}$. Then the positive roots are 
$\frac{1}{2}(\gamma_{i}-\gamma_{j})$, $1\leq i<j\leq r$, as in \ref{SS: tube type}. The coordinate ring $\CC[A]$ can be identified with the ring $\CC[x_1^{\pm 1}, \ldots, x_r^{\pm 1}]/ (\prod x_i -1)$ in such a way that 
$e^{\epsilon_{i}}=e^{\gamma_{i}/2}$ corresponds to $x_{i}$. In particular, 
we have a projection homomorphism $\pi: \CC[x_1, \ldots, x_r] \to \CC[A]$.
If we let the symmetric group $S_{r}$ act on $\CC[x_1, \ldots, x_r]$ 
as usual, then $\pi$ induces a homomorphism of rings of invariants 
$\pi: \CC[x_1, \ldots, x_r]^{S_{r}} \to \CC[A]^{W}$.

Let 
$\Lambda = \{(\lambda_1,\ldots,\lambda_r)\in\ZZ^r \mid
ÊÊÊ \lambda_1\geq \lambda_2\geq \cdots \geq \lambda_r\geq 0\}.$
If $\lambda\in \Lambda$, then $\sum_{i=1}^{r}\lambda_{i}\gamma_{i}$
is an element of $2X^{*}(A)^{+}$; thus we have a natural map 
$\Lambda\rightarrow 2X^{*}(A)^{+}$, which is surjective. 
By abuse of notation, if $\lambda\in \Lambda$, we will denote its image in  
$2X^{*}(A)^{+}$ also by $\lambda$.
With this convention, the image of the monomial $x^{\lambda}:=x_{1}^{\lambda_{1}}\cdots x_{r}^{\lambda_{r}}$ under the mapping $\pi$ is the element 
$e^{\lambda/2}$. It follows that for $\lambda\in \Lambda$,
the function $m_{\lambda}\in \CC[A]$ that was defined in \ref{p.jack} (just before the proposition) is the image of the monomial symmetric function corresponding to the partition $\lambda$.

Corresponding to any partition $\lambda\in \Lambda$ there is an element $P_{\lambda}^{( \alpha )}  \in \CC(\alpha)[x_1, \ldots, x_r]$, called a Jack polynomial
(or Jack symmetric function).  Here $\alpha$ is a parameter, not to be confused with a root (unfortunately it is customary to denote this parameter by $\alpha$).  These polynomials are defined in 
\cite[Section~10]{Macdonald}; here we recall the formulation of  \cite{KnopSahi1997}.  If $\alpha$ is such that $1/\alpha$ is a non-negative integer, then $\Delta^{1/\alpha}(x) := \prod_{i \neq j} (1 - x_i/x_j)^{1 / \alpha}$ is an element of the Laurent polynomial ring $\CC[x_1^{\pm 1}, \ldots, x_r^{\pm 1}]$.  
Define an inner product on $\CC[x_1, \ldots, x_r]$ by the rule
$$
\langle f,g \rangle_{\alpha} = [f(x) g(x^{-1}) \Delta^{1/\alpha}(x)]_0
$$
where the subscript $0$ denotes taking the constant term of a Laurent polynomial.  This inner product can be defined in an alternative way using integration; the alternative definition makes sense for all nonnegative real $\alpha$; see \cite{Macdonald}.   Then the Jack polynomials are characterized by the fact that the coefficient in $P_{\lambda}^{( \alpha )}$ of  the monomial symmetric function corresponding to $\lambda$ is $1$, and by the fact that if $\lambda \neq \mu$, then
$$
\langle P_{\lambda}^{(\alpha)}, P_{\mu}^{(\alpha)} \rangle_{\alpha} = 0
$$
for all  $\lambda, \mu\in \Lambda$ and  $\alpha$ with $1/\alpha \in \NN$ (see \cite{KnopSahi1997}).  The functions $P_{\lambda}^{(1)}$ are the Schur functions $s_{\lambda}$.   A version of the following proposition appears in \cite{Ruitenburg1989}.

\begin{prop}
Let $\lambda \in \Lambda$ and identify $\lambda$ with its image 
in $2 X^*(A)^+$.
Then under the map $\pi: \CC[x_1, \ldots, x_r]^{S_{r}}\rightarrow \CC[A]^{W}$,
the Jack polynomial $P_{\lambda}^{(2/m)} $  
is mapped to the function $P_{\lambda}$.
\end{prop}

\begin{proof}
Let $f$ and $g$ be homogeneous polynomials in $\CC[x_1, \ldots, x_r]$.  If $f$ and $g$ have the same degree, then $\langle f, g \rangle_{2/m} = \langle \pi(f), \pi(g) \rangle_A$; if $f$ and $g$ have different degrees, then $\langle f, g \rangle_{2/m} = 0$. It follows that the images of the Jack polynomials
$P_{\lambda}^{(2/m)}$ satisfy the properties of Proposition~\ref{p.jack}.  
Since these properties characterize the $P_{\lambda}$, the result follows.
\end{proof}

\subsection{}\label{SS: Stanley}
One can expand the product of Jack polynomials as a sum of Jack polynomials
$$
P_{\lambda}^{(\alpha)}  P_{\mu}^{(\alpha)} = \sum f_{\lambda \mu}^{\nu}(\alpha) P_{\nu}^{(\alpha)},
$$ 
where the $f_{\lambda \mu}^{\nu}(\alpha)$ are rational functions of $\alpha$. If $\alpha=1$, then the $f_{\lambda \mu}^{\nu}(\alpha)=c_{\lambda \mu}^{\nu}$ are Littlewood-Richardson coefficients for the group $GL_r(\CC)$; if
$\alpha=2/m$, then the $f_{\lambda \mu}^{\nu}(\alpha)=a_{\lambda \mu}^{\nu}$
are the structure coefficients for the multiplication of the spherical functions. Stanley conjectured that the $f_{\lambda \mu}^{\nu}(\alpha)$ have a certain positivity property, which we can formulate as follows. 
(Stanley's formulation used another inner product, which we do not want to define here, but the equivalence of the formulations follows from \cite[Ch.~VI, Section~10]{Macdonald}.)   For a partition $\lambda$, let 
$c_{\lambda}(\alpha)$ and $c'_{\lambda}(\alpha)$ be the polynomials in $\alpha$
with non-negative integer coefficients as defined in \cite[Ch.~VI, Section~10, (10.21)]{Macdonald}. It is known that $c_{\lambda}(\alpha) c_{\mu}(\alpha) c'_{\nu}(\alpha) f_{\lambda \mu}^{\nu}(\alpha)$ is a polynomial with integer coefficients. This follows from another conjecture of Stanley and Macdonald
(see \cite[Ch.~VI, Section~10, (10.26?) and (10.33)]{Macdonald}), which was proved by Knop and Sahi in \cite{KnopSahi1997}, and also (in part) by Lapointe and Vinet \cite{LapointeVinet1995}.

\begin{guess}[{Stanley~\cite[Conjecture 8.3]{Stanley1989}}]
For fixed partitions $\lambda,\mu,\nu\in \Lambda$, the polynomial 
$c_{\lambda}(\alpha) c_{\mu}(\alpha) c'_{\nu}(\alpha) f_{\lambda \mu}^{\nu}(\alpha)$ has non-negative integer coefficients.
\end{guess}

We will only need the following consequence of Stanley's conjecture.

\begin{cor}  
Assume that Stanley's conjecture holds.  Then 
for fixed partitions $\lambda,\mu,\nu\in \Lambda$ the following are equivalent:
\begin{itemize}
\item[(i)] $f_{\lambda \mu}^{\nu}(\alpha)$ is non-zero for all positive real values of $\alpha$.
\item[(ii)]  The Littlewood-Richardson coefficient $c_{\lambda \mu}^{\nu}$ is non-zero. 
\item[(iii)]  $a_{\lambda \mu}^{\nu}$ is non-zero.
\end{itemize}
\end{cor}

\begin{proof}
Clearly, (i) implies (ii) and (iii).  Conversely, suppose (ii) holds.
We have $f_{\lambda \mu}^{\nu}(1) = c_{\lambda \mu}^{\nu}$.  If this is non-zero, then some coefficient of the polynomial $c_{\lambda}(\alpha) c_{\mu}(\alpha) c'_{\nu}(\alpha) f_{\lambda \mu}^{\nu}(\alpha)$ must be non-zero.  Stanley's conjecture implies that this polynomial is non-zero for all positive real values of $\alpha$, proving (i).  The proof that (iii) implies (i) is similar.
\end{proof}

\subsection{}\label{SS:Main Theorem}
We now return to the notation of Section~\ref{S:hermitian}.   
By Proposition \ref{reduction}, Conjecture~B
is true iff it is true for pairs of tube type, so we assume the pair $(\g, \k)$ is of tube type.   
Recall that $\Lambda$ denotes the set of partitions of length at most $r$, identified with a subset of $\widehat{K}$ by sending $\lambda$ to  $V_\lambda$, the irreducible $K$-representation of lowest
weight $- \sum \lambda_i \gamma_i$. 
We can decompose $S = \CC[\p^+]$ into $K$-isotypic components:
$$
ÊÊÊ S  \simeqÊ \bigoplus_{\lambda \in \Lambda} S_\lambda,
$$
where $S_\lambda$ is isomorphic to  $V_\lambda$.   Note that in this case $K$ is not semisimple; it has a one-dimensional center.  However, we can replace $K$ by $K'$, the simply connected cover of the derived group $[K,K]$, and $M$ by the fixed point set $M' \subseteq K'$ of the corresponding involution.  
We have surjective group homomorphism $K'\times \CC^{\times}\rightarrow K$
and hence we may view representations of $K$ as representations of 
$K'\times \CC^{\times}$.  
Thus, we can view the above decomposition as a $K'\times \CC^{\times}$-module decomposition of $S$.

Since $M$ has finite intersection with the center of $K$, the Lie algebras of $M$ and $M'$ coincide.  Decomposing into eigenspaces of the involution $\tau$ gives
$\k' = \m \oplus \s'$ and $\k = \m \oplus \s$, and $\s = \s' \oplus \z$ where $\z$ is the center of $\k$.  
In fact, we have $\k'_{\RR} = \m_{\RR} \oplus \s'_{\RR}$, and we can take our maximal abelian
subspace of $\s'_{\RR}$ to be $\a'_{\RR} := \h_{\RR}^- \cap \s'_{\RR}$.  We obtain a corresponding 
algebraic torus $A'$ of $K'$.  By Moore's result \cite{Moore1964}  the restricted roots of $(\k',\m)$ 
are of the form $\frac{1}{2}(\gamma_i-\gamma_j)$.  The integrality conditions of
\cite[Ch.~VII \S 8]{HelgasonDG} imply that $\epsilon_j := \frac{1}{2} \gamma_j \in X^*(A')$.
Therefore, the analysis of Sections \ref{SS:root}-\ref{SS: Stanley} applies to $\CC[K'/M']$.

\begin{thm}
Suppose that Stanley's conjecture is true.  Then for 
$\lambda,\mu,\nu \in \Lambda$,
$$
S_\nu \subseteq S_{\lambda} \cdot S_{\mu}
$$
if and only if the Littlewood-Richardson coefficient $c_{\lambda \mu}^{\nu}$ is non-zero.
\end{thm}

\begin{proof}
By Proposition~\ref{reduction}, we may assume $(\g, \k)$ is of tube type.  In this case, recall from \ref{SS: tube type} that we have an inclusion
$$
S = \bigoplus_{\lambda \in \Lambda} S_\lambda \hookrightarrow \bigoplus_{\lambda \in \widehat{K}} R_\lambda,
$$
where $R = \CC[K/M]$ is a multiplicity free $K$-representation.   Replacing $R$ by $R' = \CC[K'/M'\times \CC^{\times}]=\CC[K'/M']\otimes \CC[t,t^{-1}]$, we can view $S$ as a subring of $R'$; in particular, if $\lambda \in \Lambda$, then
$S_{\lambda} = R'_{\lambda}$.  As above, for each $\lambda \in \widehat{K}$ we have an $M$-invariant element $\phi_{\lambda}$ of $R'_{\lambda}$.  If $\lambda$ and $\mu$ are in $\Lambda$, then, since $S$ is a subring of $R'$, we have
$$
\phi_{\lambda} \phi_{\mu} = \sum a_{\lambda \mu}^{\nu} \phi_{\nu}.
$$
By Ruitenburg's result, we have $S_\nu \subseteq S_{\lambda} \cdot S_{\mu}$ if and only if
$a_{\lambda \mu}^{\nu}\neq 0$.
By our discussion in the previous subsection,
$a_{\lambda \mu}^{\nu} = f_{\lambda \mu}^{\nu}(2/m)$ .  Therefore, assuming Stanley's conjecture, $a_{\lambda \mu}^{\nu}\not=0$ if and only if
$c_{\lambda \mu}^{\nu} \neq 0$.  The theorem follows.
\end{proof}

\section{Multiplicities and Littlewood-Richardson coefficients}\label{S:tensor}

\subsection{}
Let $(\g,\k)$ be an irreducible Hermitian symmetric pair of rank $r$ and
let $\{\gamma_1,\ldots,\gamma_r\}$ be the set of strongly orthogonal
roots in $\p^+$ as defined in Section~\ref{S:hermitian}.
Let $\Lambda=\{(\lambda_1,\ldots,\lambda_r) \in \ZZ^r |
\lambda_1\geq \ldots\geq \lambda_r\geq 0\}$.
In this section, if $\lambda\in \Lambda$ let $\tilde{\lambda}$ denote the
$\Delta_c^+$-dominant integral weight given by
$
Ê \tilde{\lambda}=\lambda_1\gamma_1+\cdots+\lambda_r\gamma_r.
$
Then define
\begin{itemize}
\item[]
$V_{\lambda}\ = \mbox{irreducible $\k$-module with lowest weight
$-\tilde{\lambda}$}$
\smallskip
\item[]
$F_\lambda^{(r)} = \mbox{irreducible $\gl_r$-module with lowest weight 
$-\lambda$}$
\end{itemize}
For $\lambda,\mu,\nu\in \Lambda$, let
$[V_{\lambda}\otimes V_{\mu},V_\nu]$ denote the multiplicity of $V_\nu$
in $V_{\lambda}\otimes V_{\mu}$.

\begin{prop}\label{equiv}
Let $(\g,\k)$ be an irreducible Hermitian symmetric pair of classical type
of rank $r$ and let $\lambda,\mu,\nu \in \Lambda$
such that $|\nu|=|\lambda|+|\mu|$.
Then
$$
Ê [V_{\lambda}\otimes V_{\mu},V_\nu]\not=0 \quad
Ê \mbox{if and only if}\quadÊ c^\nu_{\lambda\mu}\not =0.
$$
\end{prop}
\noindent
As an immediate consequence we obtain:
\begin{cor}
Let $(\g,\k)$ be an irreducible Hermitian symmetric pair of classical
type of rank $r$. Then Conjecture~A is true for $(\g,\k)$ if and only
if Conjecture~B is true for $(\g,\k)$.
\end{cor}

\subsection{}
Before we prove Proposition~\ref{equiv} we recall some known facts
about Littlewood-Richardson coefficients that are needed in the proof.
We first need a little bit more notation about partitions.
We use the same notation as in Macdonald's book \cite{Macdonald}.
A partition is a sequence
$\lambda=(\lambda_1,\lambda_2,\ldots,\lambda_r,\ldots)$
of non-negative integers such that $\lambda_1\geq \lambda_2\geq \cdots
\geq \lambda_r \geq \cdots$ and $\lambda_i=0$ for $i\gg 1$. The integer
$\ell(\lambda)=\min\{i |\lambda_i\not=0\}$ is the length of $\lambda$
and $|\lambda|=\sum_i\lambda_i$ is the size of $\lambda$. The partition
conjugate to $\lambda$ is the partition $\lambda'$ given by
$(\lambda')_i=\#\{\lambda_j \mid \lambda_j\geq i\}$.
Let $\lambda,\mu,\nu$ be partitions of length at most $r$. By the
Littlewood-Richardson rule it follows that for any $n\geq r$,
$$
Ê [F_\lambda^{(n)} \otimes F_\mu^{(n)},F_\nu^{(n)}]=
Ê [F_\lambda^{(r)} \otimes F_\mu^{(r)},F_\nu^{(r)}]=c^{\nu}_{\lambda\mu}.
$$
Another interpretation of Littlewood-Richardson coefficients is in terms of Schur  functions. For partition $\lambda$ and $\mu$, we have 
$s_{\lambda}s_{\mu}=\sum_{\nu}c^{\nu}_{\lambda\mu}s_{\nu}$.
A simple, but very useful necessary condition for $c^\nu_{\lambda\mu}\not=0$
is that $|\lambda|+|\mu|=|\nu|$.
We will also need that
$c^\nu_{\lambda\mu}=c^{\nu'}_{\lambda'\mu'}$.
This follows since there is an involution $\omega$ on the ring of symmetric functions
such that $\omega(s_\lambda)=s_{\lambda'}$. Finally, we will need that 
$$
   \mbox{$c^{\nu}_{\lambda\mu}\not=0$}\quad   
   \mbox{if and only if}
  \quad
Ê \mbox{$c^{2\nu}_{2\lambda, 2\mu}\not=0$.}
$$
This is a special case of Klyachko's saturation conjecture 
\cite{Klyachko1998}, which is now a theorem of Knutson and Tao \cite{KnutsonTao1999} and also of Derksen and Weyman~\cite{DerksenWeyman1990}.

\subsection{}{\it Proof of Proposition~\ref{equiv}.}
We will prove the proposition case by case.

\smallskip
\noindent
{\bf Case} $\g_\RR=\su(p,q)$.
It is slightly more convenient to work withÊ $\g_\RR=\u(p,q)$ instead
of $\su(p,q)$.
In this case, $(\g,\k)=(\gl_n,\gl_p\times \gl_q)$ and $r=p$. Here we assume
that $p\leq q$.
From Table~\ref{gammas}, if $\lambda=(\lambda_1,\ldots,\lambda_p)$
then
$
\tilde{\lambda}=(\lambda_1,\ldots,\lambda_p,0,\ldots,0,-\lambda_p,\ldots,-\lambda_1).
$
Hence $V_\lambda\simeq F_\lambda^{(p)} \boxtimes \big(F_\lambda^{(q)}\big)^*$
and it follows that $[V_\lambda\otimes V_\mu, V_\nu]=(c^\nu_{\lambda\mu})^2.$
The equivalence is now obvious.

\medskip
\noindent
{\bf Case $\g_\RR=\sp(n,\RR)$.}
In this case, $(\g,\k)=(\sp_{2n},\gl_n)$ and $r=n$.
From Table~\ref{gammas}, if $\lambda=(\lambda_1,\ldots,\lambda_n)$
then $\tilde{\lambda}=(2\lambda_1,\ldots,2\lambda_n)=2\lambda.$
Hence $V_\lambda \simeq F_{2\lambda}^{(n)}$ and it follows that
$[V_\lambda\otimes V_\mu, V_\nu]=c^{2\nu}_{2\lambda,2\mu}.$
The equivalence is now an immediate consequence of the saturation conjecture.

\medskip
\noindent
{\bf Case $\g_\RR=\so^*(2n)$}.
In this case, $(\g,\k)=(\so_{2n},\gl_n)$ and $r=\lfloor n/2\rfloor$.
From Table~\ref{gammas}, if $\lambda=(\lambda_1,\ldots,\lambda_r)$ then
$
\tilde{\lambda}=(\lambda_1,\lambda_1,\lambda_2,\lambda_2,\ldots)=(2\lambda')'.
$
Hence $V_\lambda \simeq F_{(2\lambda')'}^{(n)}$ and it follows that
$
[V_\lambda\otimes V_\mu, V_\nu]=c^{(2\nu')'}_{(2\lambda')',(2\mu')'}=c^{2\nu'}_{2\lambda',2\mu'}.
$
The equivalence is now a consequence of the saturation conjecture and
the fact that $c^{\nu'}_{\lambda'\mu'}=c^\nu_{\lambda\mu}$.

\medskip
\noindent
{\bf Case $\g_\RR=\so(n,2)$}.
In this case, $(\g,\k)=(\so_{n+2},\so_n\oplus\CC)$ and $r=2$. We may assume
that $n\geq 4$ since $\so(3,2)\simeq \sp(2,\RR)$.
Let $\{\epsilon_1,\ldots,\epsilon_{\lfloor n/2\rfloor}\}$
be the canonical basis for the standard Cartan subalgebra of $\so_n$.
Extend this basis to the standard basis (with non-standard labeling)
$\{\epsilon_0,\epsilon_1,\ldots,\epsilon_{\lfloor n/2\rfloor}\}$ of
$\so_{n+2}$. With this convention, $\gamma_1=\epsilon_0+\epsilon_1$
and $\gamma_2=\epsilon_0-\epsilon_1$. So,
if $\lambda=(\lambda_1,\lambda_2)$ then
$\tilde{\lambda}=(\lambda_1+\lambda_2)\epsilon_0+
(\lambda_1-\lambda_2)\epsilon_1$.
Thus, as a $\k=\so_{n}\oplus \CC$-module,
$V_\lambda \simeq E_{(\lambda_1-\lambda_2)\epsilon_1}^{(n)}\boxtimes \CC_{-\lambda_1-\lambda_2}$,
where $E_{(\lambda_1-\lambda_2)\epsilon_1}^{(n)}$ is the irreducible
$\so_n$-module with lowest weight $-(\lambda_1-\lambda_2)\epsilon_1$.
In the following we will show that $[V_\lambda\otimes V_\mu, V_\nu]=c^\nu_{\lambda\mu}$.

\begin{lemma}
Suppose that $n\geq 4$.
If $a$ and $b$ are non-negative integers then
$$
E_{a\epsilon_1}^{(n)}\otimes E_{b\epsilon_1}^{(n)}\ \simeq\
\bigoplus_{k=0}^{b}\ \bigoplus_{l=0}^{b-k}
ÊE_{(a+b-2k-l)\epsilon_1+l\epsilon_2}^{(n)}\ .
$$
\end{lemma}

\noindent{\it Proof of the lemma\/.} The assumption $n\geq 4$ makes
sure that we are in the ``stable range''. We can then use the branching
rule \cite[Theorem 2.1.2]{HoweTanWillenbring2004} to decompose
the tensor product $E_{a\epsilon_1}^{(n)}\otimes E_{b\epsilon_1}^{(n)}$. 
The lemma then follows by the usual Pieri rule for Littlewood-Richardson coefficients.\hfill\qed

\bigskip
\noindent
We now return to the proof of Proposition~\ref{equiv}.
Let $\lambda,\mu,\nu\in\Lambda$ such that $|\lambda|+|\mu|=|\nu|$.
We note that $F_\lambda^{(2)}\simeqÊ F_{(\lambda_1-\lambda_2,0)}^{(2)}
\otimes F_{(\lambda_2,\lambda_2)}^{(2)}$ and
$F_{(\lambda_2,\lambda_2)}^{(2)}$ is one-dimensional,
which implies that we can calculate $c^\nu_{\lambda\mu}$
by using the Pieri rule.
If we put $a=\lambda_1-\lambda_2$, $b=\mu_1-\mu_2$ and $c=\nu_1-\nu_2$
then
$$
c^\nu_{\lambda\mu}=
\begin{cases}
1 & \mbox{if $c=a+b-2k$ for some non-negative integer $k$}\\
0 & \mbox{otherwise}
\end{cases}
$$
By the lemma above,
$[E_{a\epsilon_1}^{(n)}\otimes E_{b\epsilon_1}^{(n)},E_{c\epsilon_1}^{(n)}]
=c^{\nu}_{\lambda\mu}$ and hence
$[V_\lambda\otimes V_\mu, V_\nu]=c^\nu_{\lambda\mu}$.
This completes the proof of Proposition~\ref{equiv}. \hfill\qed

\subsection{}
We conjecture that Proposition~\ref{equiv} also holds in the two exceptional cases. 

\medskip
\noindent
{\bf Case $\g_\RR= E\, {I\!I\!I}$}.
In this case, $(\g,\k)=(\frak e_6,\so_{10}\oplus \CC)$ and $r=2$.
Let $\omega_1,\ldots,\omega_6$ be the fundamental weights of $\frak e_6$ (Bourbaki ordering).
By Table~2, $\gamma_1=\omega_2$, $\gamma_2=\omega_1-\omega_2+\omega_6$.
Thus if $\lambda=(\lambda_1,\lambda_2)$ then
$\tilde{\lambda}=\lambda_2\omega_1+(\lambda_1-\lambda_2)\omega_2+\lambda_2\omega_6$. Let $\varpi_1,\dots,\varpi_5$ be the fundamental weights of $\so_{10}$ (Bourbaki ordering). If $\lambda=(\lambda_1,\lambda_2)\in \Lambda$ then as a
$\k=\so_{10}\oplus \CC$-module, $V_\lambda \simeq E_{(\lambda_1-\lambda_2)\varpi_1+\lambda_2\varpi_5}\boxtimes \CC_{-\lambda_1-\lambda_2}$, where $E_{(\lambda_1-\lambda_2)\varpi_1+\lambda_2\varpi_5}$ denotes the irreducible $\so_{10}$-module of lowest weight
$-(\lambda_1-\lambda_2)\varpi_1-\lambda_2\varpi_5$. We conjecture that if 
$\lambda,\mu,\nu\in\Lambda$ with $|\lambda|+|\mu|=|\nu|$,
then $[V_\lambda\otimes V_\mu, V_\nu] = c^{\nu}_{\lambda\mu}$.

\medskip
\noindent
{\bf Case $\g_\RR= E\, {V\!I\!I}$}.
In this case, $(\g,\k)=(\frak e_7,\frak e_6 \oplus \CC)$ and $r=3$.
Let $\omega_1, \ldots, \omega_7$ be the fundamental weights of $\frak e_7$ (Bourbaki ordering).
By Table~2, $\gamma_1=\omega_1$, $\gamma_2=-\omega_1+\omega_6$ and 
$\gamma_3=-\omega_6+2\omega_7$. If $\lambda=(\lambda_1,\lambda_2,\lambda_3)$ then $\tilde{\lambda}=(\lambda_1-\lambda_2)\omega_1+(\lambda_2-\lambda_3)\omega_6+2\lambda_3\omega_7$. Let $\varpi_1,\ldots,\varpi_6$ be the 
fundamental weights of $\frak e_6$ (Bourbaki ordering). If $\lambda=(\lambda_1,\lambda_2,\lambda_3)\in \Lambda$ 
then as a $\k=\frak e_{6}\oplus \CC$-module,
$V_\lambda \simeq W_{(\lambda_1-\lambda_2)\varpi_1+
(\lambda_2-\lambda_3)\varpi_6} \boxtimes \CC_{-\lambda_1-\lambda_2-\lambda_3}$,
where  $W_{(\lambda_1-\lambda_2)\varpi_1+
(\lambda_2-\lambda_3)\varpi_6}$ denotes the irreducible $\frak e_6$-module 
with lowest weight  $-(\lambda_1-\lambda_2)\varpi_1-
(\lambda_2-\lambda_3)\varpi_6$.
We conjecture that if $\lambda,\mu,\nu\in\Lambda$ with $|\lambda|+|\mu|=|\nu|$,
then 
$$[V_\lambda\otimes V_\mu, V_\nu]=
\frac{c^\nu_{\lambda\mu}(c^\nu_{\lambda\mu}+1)}{2}.
$$
We arrived at this conjecture by computing examples using the computer algebra package LiE \cite{LiE1992}.

\begin{table}[h]
\begin{center}
\begin{tabular}{|l|l|l|}
\hline
$K$ & $X=\p^{+}$ & $ [V_{\lambda}\otimes V_{\mu},V_{\nu}]$\\ 
\hline
$GL_{p}(\CC)\times GL_{q}(\CC)$ & $\CC^{p}\otimes (\CC^{q})^{*}$ & $(c_{\lambda\mu}^{\nu})^{2}$ \\[1.0ex]
$GL_{n}(\CC)$ & $S^{2}(\CC^{n})$ & $c_{2\lambda,2\mu}^{2\nu}$ \\[1.0ex]
$GL_n(\CC)$ & $\wedge^2(\CC^n)$ & $c_{(2\lambda')',(2\mu')'}^{(2\nu')'}=c_{2\lambda,2\mu}^{2\nu}$\\[1.0ex]
$SO_{n}(\CC)\times \CC^{\times}$ & $ \CC^n$ & $c_{\lambda\mu}^{\nu}$ \\[1.0ex]
$Spin_{10}(\CC) \times \CC^\times$  & $\CC^{16}$ (spin) &$c_{\lambda\mu}^{\nu}$
(?)\\[1.0ex] 
$E_{6}(\CC)\times \CC^\times$ & $\CC^{27}$ (min)& $\frac{c_{\lambda\mu}^{\nu}(c_{\lambda\mu}^{\nu}+1)}{2}$ (?)\\[1.0ex]
\hline
\end{tabular}
\end{center}
\vspace{.5pc}
\caption{Multiplicities $[V_{\lambda}\otimes V_{\mu},V_{\nu}]$
in terms of Littlewood-Richardson coefficients}\label{T:multiplicities}
\end{table}

\end{document}